\documentclass{article}

\usepackage{graphicx}

\usepackage{amssymb}

\begin{document}

\title{ New properties of the Lerch's transcendent}

\author{\ \\
E.M. Ferreira,  A.K. Kohara, \\  \  \\
{\em Instituto de F\'{\i}sica, Universidade Federal do Rio de Janeiro,} \\
{\em 21941-972, Rio de Janeiro, Brasil}\\ \ \\
and \\ \ \\
J. Sesma\thanks{Corresponding author. Email: javier@unizar.es} \\ \ \\
{\em Departamento de F\'{\i}sica Te\'{o}rica, Facultad de Ciencias,} \\
{\em 50009, Zaragoza, Spain}}

\maketitle

\begin{abstract}
A new representation of the Lerch's transcendent $\Phi(z,s,a)$, valid for positive integer $s=n=1, 2, \ldots$ and for $z$ and $a$ belonging to certain regions of the complex plane, is presented. It allows to write an equation relating $\Phi(z,n,a)$ and $\Phi(1/z,n,1-a)$, which provides an expansion of $\Phi(z,n,a)$ as a power series of $1/z$, convergent for $|z|>1$.
\end{abstract}

{\bf AMS Subject Clasification:} 11M35; 30E15

\section{Introduction}

The Lerch's transcendent $\Phi(z,s,a)$, also known as Hurwitz-Lerch zeta function, is defined by its series representation \cite[Sec.~1.11, Eq.~(1)]{erde}
\cite[Eq.~25.14.1]{dlmf}
\begin{equation}
\Phi(z,s,a) = \sum_{m=0}^\infty\frac{z^m}{(a+m)^s}\,,   \label{i1}
\end{equation}
provided
\begin{equation}
a\neq 0, -1, -2, \ldots; \qquad |z|<1; \qquad |z|=1, \quad \Re s>1\,.   \label{i2}
\end{equation}
The restriction on the values of $a$ guarantees that all terms of the series in the right-hand side are finite. Obviously,
the series is convergent if $|z|<1$, independently of the value of $s$, or if $|z|=1$ and $\Re s>1$.
For other values of its arguments, $\Phi(z,s,a)$ is defined by analytic continuation. This is achieved by
means of integral representations, the most common of them being \cite[Eq.~25.14.5]{dlmf}
\begin{equation}
\Phi(z,s,a) = \frac{1}{\Gamma(s)}\int_0^{\infty}\frac{t^{s-1}\,{e}^{-at}}{1-z\,{e}^{-t}}\,{d}t \,,  \label{i3}
\end{equation}
whenever
\begin{equation}
\Re s>0\,, \qquad \Re a>0\,, \qquad z\in \mathbb{C}\setminus [1,\infty)\,.  \label{i4}
\end{equation}
The conditions imposed on $a$ and $z$ ensure the regularity of the integrand in the right-hand side of (\ref{i3}).
The restriction on $s$ allows to prove the equivalence of the representations (\ref{i1}) and (\ref{i3}) in their
common region of validity (see, for instance, \cite[Lemmas 2.1 and 2.2]{guil}).
A thorough discussion of the analytic continuation of $\Phi$, as a multivalued function of three complex
variables, and of its ``singular strata" can be found in a recent paper by Lagarias and Li \cite{laga}, where the monodromy functions
describing the multivaluedness are computed. (Be aware that, in the notation used in Ref.~\cite{laga},
the two first arguments of $\Phi$ are transposed, as compared with the notation used in Refs.~\cite{erde}
and \cite{dlmf}.)

In the course of a research on the use of dispersion relations in the study of elementary particles \cite{fer1,fer2}, we have
encountered what we believe to be a new representation of the Lerch's transcendent $\Phi(z,s,a)$ for positive integer
values of the second argument, $s=n=1, 2, \ldots$. This is our first result, presented as Theorem 1 in Sec.~2.
Such a representation allows to unveil, as a second result reported in Sec.~3 as Theorem 2,
a property of $\Phi$ not noticed before. This property, in turn, provides our third result, expressed in Corollary 1 of Sec.~4,
consisting of an expansion of $\Phi(z,n,a)$ in powers of $1/z$, convergent for $|z|>1$. The proofs of the three results, followed by pertinent remarks,  are presented in Secs.~2, 3 and 4, respectively. Some comments are added in Sec.~5.

\section{A new representation}

Let $\mathbf{D}$ denote the open unit disc in the complex plane, cut along the negative real semiaxis, that is,
\begin{equation}
z\in\mathbf{D} \qquad \Rightarrow \qquad z\in \mathbb{C}, \quad 0<|z|<1, \quad -\pi<\arg(z)<\pi.   \label{ii1}
\end{equation}
A new representation of the Lerch's transcendent $\Phi(z,s,a)$, for positive integer values of its second argument, is provided by the following

\medskip
\noindent {\bf Theorem 1.} {\em Let us assume $z\in\mathbf{D}$ and denote
\begin{equation}
\varphi=\arg(-\ln z).   \label{ii2}
\end{equation}
For positive integer $n=1, 2, \ldots$, and complex $a$ such that $\Re[(a-1){e}^{{i}\varphi}]<0$,
the Lerch's transcendent admits the representation
\begin{equation}
\Phi(z,n,a) = \frac{(-1)^{n-1}}{(n-1)!}\left\{P\int_0^{\infty {e}^{{i}\varphi}}\frac{t^{n-1}\,{e}^{at}}{z\,{e}^t-1}\,{d}t
+ \pi\,\frac{\partial^{n-1}}{\partial a^{n-1}}\left(z^{-a}\,\cot(\pi a)\right)\right\},   \label{ii3}
\end{equation}
where the symbol $P$ stands for the Cauchy principal value of the path integral along the ray $\arg(t)=\varphi$.
}

\noindent {\bf Proof.} Our starting point is an integral encountered as we were writing dispersion relations for proton-proton and proton-antiproton scattering at high energies \cite{fer1,fer2}. According to the definition of principal value of an integral,
\begin{eqnarray}
P\int_0^{\infty {e}^{{i}\varphi}}\frac{t^{n-1}\,{e}^{at}}{z\,{e}^{t}-1}\,{d}t & = & \lim_{\varepsilon\to 0^+}
\Bigg\{-\int_0^{-\ln z-\varepsilon {e}^{{i}\varphi}}\frac{t^{n-1}\,{e}^{at}}{1-z\,{e}^{t}}\,{d}t  \nonumber  \\
 & & \hspace{30pt} +\,z^{-1}\int_{-\ln z+\varepsilon {e}^{{i}\varphi}}^{\infty {e}^{{i}\varphi}}\frac{t^{n-1}\,
 {e}^{(a-1)t}}{1-z^{-1}\,{e}^{-t}}\,{d}t\Bigg\}\,.  \label{iii1}
\end{eqnarray}
By replacing, in the right-hand side, $1/(1-z\,{e}^{t})$ and $1/(1-z^{-1}\,{e}^{-t})$ by their respective
geometric series expansions we obtain
\begin{eqnarray}
P\int_0^{\infty {e}^{{i}\varphi}}\frac{t^{n-1}\,{e}^{at}}{z\,{e}^{t}-1}\,{d}t & = & \sum_{m=0}^\infty\Bigg[
\lim_{\varepsilon\to 0^+}\Bigg\{-z^m\int_0^{-\ln z-\varepsilon {e}^{{i}\varphi}}{e}^{(m+1/2)t}\,G(t)\,{d}t  \nonumber  \\
& & \hspace{10pt} +\,z^{-m-1}\int_{-\ln z+\varepsilon {e}^{{i}\varphi}}^{\infty {e}^{{i}\varphi}} {e}^{-(m+1/2)t}\,
G(t)\,{d}t\Bigg\}\Bigg]\,,  \label{iii3}
\end{eqnarray}
with the notation
\begin{equation}
G(t) \equiv t^{n-1}\,{e}^{(a-1/2)t}\,.   \label{iii4}
\end{equation}
Repeated integration by parts gives
\begin{eqnarray}
P\int_0^{\infty {e}^{{i}\varphi}}\frac{t^{n-1}\,{e}^{at}}{z\,{e}^{t}-1}\,{d}t & = & \nonumber  \\
& & \hspace{-50pt} \sum_{m=0}^\infty\Bigg[\lim_{\varepsilon\to 0^+} \Bigg\{ -z^m\sum_{k=0}^\infty
\left[\frac{(-1)^k{e}^{(m+1/2)t}}{(m+1/2)^{k+1}}\,
\frac{{d}^kG(t)}{{d}t^k}\right]_{t=0}^{t=-\ln z-\varepsilon {e}^{{i}\varphi}}  \nonumber  \\
&  &  -\,z^{-m-1}\sum_{k=0}^\infty\left[\frac{{e}^{-(m+1/2)t}}{(m+1/2)^{k+1}}\,
\frac{{d}^kG(t)}{{d}t^k}\right]_{t=-\ln z+\varepsilon {e}^{{i}\varphi}}^{t\to\infty {e}^{{i}\varphi}}\Bigg\}\Bigg] \nonumber \\
& = & \sum_{m=0}^\infty z^m\sum_{k=0}^\infty\frac{(-1)^k}{(m+1/2)^{k+1}}\,\left.\frac{{d}^kG(t)}{{d}t^k}\right|_{t=0}  \nonumber  \\
& & \hspace{-10pt} +\,z^{-1/2}\sum_{m=0}^\infty\sum_{l=0}^\infty\frac{2}{(m+1/2)^{2l+2}}\,\left.
\frac{{d}^{2l+1}G(t)}{{d}t^{2l+1}}\right|_{t=-\ln z},   \label{iii5}
\end{eqnarray}
where use has been made of the condition $\Re[(a-1){e}^{{i}\varphi}]<0$. Obviously,
\begin{eqnarray}
\left.\frac{{d}^kG(t)}{{d}t^k}\right|_{t=0} & = & \left\{ \begin {array}{ll} 0\,, & {\rm if} \; k<n-1\,, \\
(k!/(k\! -\! n\! +\! 1)!)\,(a\! -\! 1/2)^{k-n+1}\,, & {\rm if} \; k\geq n-1\,,\end{array}\right. \label{bi14} \\
\left.\frac{{d}^{k}G(t)}{{d}t^{k}}\right|_{t=-\ln z} & = & \sum_{j=0}^{{\rm min}\{k, n-1\}}\Bigg[{k \choose j}
\frac{(n\! -\! 1)!}{(n\! -\! 1\! -\! j)!}\,(-\ln z)^{n-1-j}  \nonumber \\
& & \hspace{100pt}(a\! -\! 1/2)^{k-j}\,z^{-a+1/2}\Bigg],     \label{iii6}
\end{eqnarray}
which, introduced in (\ref{iii5}), give
\begin{eqnarray}
P\int_0^{\infty {e}^{{i}\varphi}}\frac{t^{n-1}\,{e}^{at}}{z\,{e}^{t}-1}\,{d}t
& = & \sum_{m=0}^\infty z^m\sum_{k=n-1}^\infty\frac{(-1)^k\,k!\,(a-1/2)^{k-n+1}}{(k-n+1)!\,(m+1/2)^{k+1}}
\nonumber  \\
& & \hspace{-100pt} +\,z^{-a}\sum_{l=0}^\infty\Bigg(\sum_{j=0}^{{\rm min} \{2l+1, n-1\}} {2l+1 \choose j}\,
\frac{(n\! -\! 1)!}{(n\! -\! 1\! -\! j)!}\,(-\ln z)^{n-1-j} \,(a\! -\! 1/2)^{2l+1-j}  \nonumber \\
& & \hspace{100pt}\sum_{m=0}^\infty\frac{2}{(m+1/2)^{2l+2}}\Bigg).   \label{iii7}
\end{eqnarray}
The last sum can be written in terms of Bernoulli numbers in the form \cite[Vol.1, Eq. (5.1.4.1)]{prud}
\begin{equation}
\sum_{m=0}^\infty\frac{2}{(m+1/2)^{2l+2}} = (2\pi)^{2l+2}\,\frac{2^{2l+2}-1}{(2l+2)!}\,\left|B_{2l+2}\right|.
\label{iii8}
\end{equation}
By recalling the series expansion of the tangent function \cite[Eq. 4.19.3]{dlmf}
\begin{equation}
\tan \alpha = \sum_{n=1}^\infty \frac{2^{2n}\,(2^{2n}-1)\,\left|B_{2n}\right|}{(2n)!}\,\alpha^{2n-1}\,,
\qquad |\alpha|<\pi/2\,,  \label{iii9}
\end{equation}
and of its derivatives
\begin{equation}
\frac{{d}^j}{{d}\alpha^j}\tan \alpha = \sum_{n=1+[j/2]}^\infty \frac{2^{2n}\,(2^{2n}-1)\,
\left|B_{2n}\right|}{(2n)!}\,\frac{(2n-1)!}{(2n-1-j)!}\,\alpha^{2n-1-j}\,, \quad |\alpha|<\pi/2\,,  \label{iii10}
\end{equation}
equation (\ref{iii7}) becomes
\begin{eqnarray}
P\int_0^{\infty {e}^{{i}\varphi}}\frac{t^{n-1}\,{e}^{at}}{z\,{e}^{t}-1}\,{d}t
& = & (-1)^{n-1}(n-1)!\sum_{m=0}^\infty \frac{z^m}{(m+a)^{n}}  \nonumber  \\
& &  \hspace{-20pt} +\,\pi\,z^{-a}\sum_{j=0}^{n-1}{n-1 \choose j}\,(-\ln z)^{n-1-j} \frac{{d}^j}{{d}a^j}
\tan(\pi(a-1/2)) \nonumber  \\
&\hspace{-30pt} = &  \hspace{-20pt}(-1)^{n-1}(n\! -\! 1)!\,\Phi(z,n,a)  - \pi\,\frac{\partial^{n\! -\! 1}}{\partial a^{n-1}}\left(z^{-a}\cot (\pi a)\right),  \label{iii11}
\end{eqnarray}
provided $|a-1/2|<1/2$, a restriction which may be relaxed, as explained in Remark 1 below. Isolation of the Lerch's transcendent completes the proof. $\square$

\medskip
\noindent {\bf Remark 1.} The condition $|a-1/2|<1/2$, necessary for the validity of (\ref{iii11}), restricts the possible complex values
of $a$ to an open disc of radius $1/2$ centered at 1/2.  Nevertheless, the representation can be continued
analytically to the complex half-plane $\Re[(a-1){e}^{{i}\varphi}]<0$, the points $a=0, -1, -2, \ldots$ being
excluded. At these points $\Phi(z,n,a)$, which for fixed $z$ is a meromorphic function of $a$, presents poles of order $n$,
as proved by Lagarias and Li \cite[Sec.~7, Theorem 7.1]{laga}. The trigonometric expression in the right-hand side of (\ref{ii3})
shows these singularities.

\section{A ``symmetry" property}

In the computation of the principal value of the integral in (\ref{iii1}), before taking the limit $\varepsilon\to 0^+$, we can
add and subtract the contribution of completing the integration path, which goes from 0 to $-\ln z-\varepsilon {e}^{{i}\varphi}$
and from $-\ln z+\varepsilon {e}^{{i}\varphi}$ to $\infty {e}^{{i}\varphi}$, with a semi-circumference (sc) of radius $\varepsilon$,
centered at $t=-\ln z$, and going from $-\ln z-\varepsilon {e}^{{i}\varphi}$ to $-\ln z+\varepsilon {e}^{{i}\varphi}$
in such a way that the point $-\ln z$ lies out of the region delimited by the resulting path, that we denote by $\Gamma$ (see Figure 1),
and the positive real semi-axis. Obviously,
\begin{equation}
\lim_{\varepsilon\to 0^+} \int_{\rm sc}\frac{t^{n-1}\,{e}^{at}}{z\,{e}^t-1}\,{d}t = {\rm sgn}(\varphi)\,
{i}\,\pi\,(-\ln z)^{n-1}\,z^{-a}.  \label{iv1}
\end{equation}
\begin{figure}
\begin{center}
\resizebox{10cm}{!}{\includegraphics{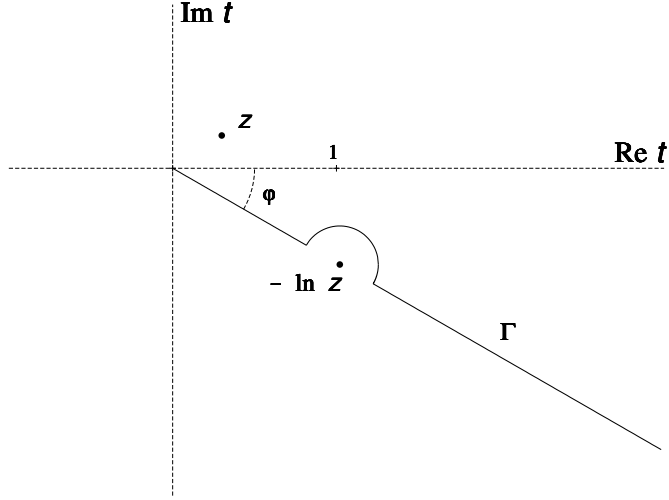}}
\caption{The integration path of the integral in the right-hand side of Eq.~(\ref{iv2}).}
\end{center}
\end{figure}
Then
\begin{equation}
P\int_0^{\infty {e}^{{i}\varphi}}\frac{t^{n-1}\,{e}^{at}}{z\,{e}^{t}-1}\,{d}t =
\int_\Gamma \frac{t^{n-1}\,{e}^{at}}{z\,{e}^{t}-1}\,{d}t-{\rm sgn}(\varphi)\, {i}\,\pi\,(-\ln z)^{n-1}\,z^{-a}.  \label{iv2}
\end{equation}
In the case of being  $\Re[(a-1){e}^{{i}c\varphi}]<0$ for all $c\in[0,1]$, the path $\Gamma$ can be deformed to coincide
with the positive real semi-axis. We obtain, in this way,
\begin{equation}
P\int_0^{\infty {e}^{{i}\varphi}}\frac{t^{n-1}\,{e}^{at}}{z\,{e}^{t}-1}\,{d}t =  z^{-1}\int_0^\infty
\frac{t^{n-1}\,{e}^{-(1-a)t}}{1-z^{-1}\,{e}^{-t}}\,{d}t - {\rm sgn}(\varphi)\, i\,\pi\,(-\ln z)^{n-1}\,z^{-a}.  \label{iv3}
\end{equation}
We have now all necessary information to prove

\medskip
\noindent {\bf Theorem 2.}
{\em For $z\in\mathbf{D}\setminus(0,1)$, $a\notin\mathbb{Z}$, and $\varphi$ as defined in (\ref{ii2}), the Lerch's transcendent with positive integer second argument presents the ``symmetry" property
\begin{eqnarray}
\Phi(z,n,a) + (-1)^{n}\,z^{-1}\,\Phi(z^{-1},n,1-a) & = &   \nonumber  \\
 &  & \hspace{-100pt} \frac{\pi\,(-1)^{n-1}}{(n-1)!}\,\frac{\partial^{n-1}}{\partial a^{n-1}}
 \left(z^{-a}\left(\cot(\pi a)-{\rm sgn}(\varphi)\,{i}\right)\right).   \label{ii4}
\end{eqnarray}
} 

\noindent{\bf Proof.} Let us consider Eq. (\ref{iv3}) above. Since the integral in the right-hand side is, up to a factor, the well known representation (\ref{i3}) of the Lerch's transcendent, we can write
\begin{equation}
P\int_0^{\infty {e}^{{i}\varphi}}\frac{t^{n-1}\,{e}^{at}}{z\,{e}^{t}-1\,}\,{d}t =  z^{-1}\,(n-1)!\,\Phi(z^{-1},n,1-a)
- {\rm sgn}(\varphi)\, {i}\,\pi\,(-\ln z)^{n-1}\,z^{-a},  \label{iv4}
\end{equation}
provided
\[
z^{-1}\notin [1, \infty) \qquad {\rm or, \,equivalently,} \qquad z\notin [0, 1]\,.
\]
Substitution of (\ref{iv4}) in the right-hand side of (\ref{ii3}), bearing in  mind that
\[
(-\ln z)^{n-1}\,z^{-a}=\frac{\partial^{n-1}}{\partial a^{n-1}}\,z^{-a}\,,
\]
gives (\ref{ii4}). The condition $\Re[(a-1){e}^{{i}c\varphi}]<0$ for all $c\in[0,1]$, required for the validity of (\ref{iv3}),
is not necessary for the property (\ref{ii4}), which can be continued analytically to the whole complex $a$-plane, the integers being excepted.
At these points, $a=\ldots, -2, -1, 0, 1, 2, \ldots$, both left and right hand sides of (\ref{ii4}) present poles of order $n$. The theorem becomes proved in this way. $\square$

\medskip
\noindent {\bf Remark 2.} It is immediate to verify that Eq. (\ref{ii4}) remains invariant under the simultaneous replacements of $z$ by $z^{-1}$
(and, consequently, ${\rm sgn}(\varphi)$ by $- {\rm sgn}(\varphi)$) and $a$ by $(1-a)$. This means that,
although it has been deduced on the assumption that $|z|<1$, the property (\ref{ii4}) applies, as it is, to the case of $|z|>1$,
provided $z\notin (1,\infty)$, and is also valid, by continuation, for $|z|=1$. The particular value $z=1$ deserves special consideration,
as it belongs to the singular stratum $\{0, 1, \infty\}$ of the double specialization of $\Phi(z,s,a)$ for positive integer $s=n$
and fixed complex $a$ \cite[Sec.~8]{laga}. In the case of $n=1$, $\Phi(1,1,a)$ and $\Phi(1,1,1-a)$ are singular and the left-hand side of (\ref{ii4})
turns out to be the difference of two infinities. For $n=2, 3, \ldots$, the fact that the Lerch's transcendent becomes the Hurwitz zeta function
when $z=1$ \cite[Eq. 25.14.2]{dlmf},
\begin{equation}
\zeta(s,a)=\Phi(1,s,a)\,, \qquad \Re s>1\,, \quad a\neq 0, -1, -2, \ldots \,,    \label{vi1}
\end{equation}
may be used to write (\ref{ii4}) in the form
\begin{equation}
\zeta(n,a)+(-1)^n\,\zeta(n,1-a)=\frac{(-1)^{n-1}\,\pi}{(n-1)!}\,\frac{{d}^{n-1}}{{d}a^{n-1}}
 \cot(\pi a)\,, \qquad n=2,3,\ldots\,,   \label{vi2}
\end{equation}
that is a known reflection property \cite[Sec. 4]{cvij}\footnote{Notice a typo in the citation in the first paragraph of Sec. 4. Instead of Ref. 6,
one should read Ref. 8.}, or in terms of polygamma functions \cite[Eq. 5.15.6]{dlmf},
\begin{equation}
\psi^{(m)}(a)-(-1)^m\,\psi^{(m)}(1-a)=-\,\pi\,\frac{{d}^{m}}{{d}a^{m}}\cot(\pi a)\,,  \qquad m=0, 1, 2, \ldots\,,  \label{vi3}
\end{equation}
a relation which follows trivially from the reflection property of the Gamma function \cite[Eq. 5.5.3]{dlmf}
\begin{equation}
\Gamma(a)\,\Gamma(1-a)=\pi/\sin(\pi a)\,.  \label{vi4}
\end{equation}

\section{A convergent expansion}

In what follows, to avoid misunderstanding, we make use of the notation
\begin{equation}
w\equiv z^{-1},  \qquad  b\equiv 1-a \,.   \label{extra1}
\end{equation}
Consequently, we have for $\varphi$, defined in (\ref{ii2}),
\begin{equation}
\varphi=\arg(\ln w).  \label{extra2}
\end{equation}
An immediate consequence of our Theorem 2 is the following

\medskip
\noindent {\bf Corollary 1.}
{\em For
\begin{equation}
|w|>1, \qquad w\notin (-\infty, 1)\cup(1, \infty),    \label{ii6}
\end{equation}
and positive integer $n$, the Lerch's transcendent $\Phi(w,n,b)$ admits the convergent expansion
\begin{equation}
\Phi(w,n,b) = \frac{\pi}{(n-1)!}\left[\frac{\partial^{n-1}}{\partial t^{n-1}}
\Big(w^{t}\big({\rm sgn}(\varphi)\,{i}-\cot(\pi t)\big)\Big)\right]_ {t=-b}
 -\,\sum_{m=1}^\infty\frac{w^{-m}}{(b-m)^n}\,. \label{ii5}
\end{equation}
}

\noindent {\bf Proof.} With the notation introduced in (\ref{extra1}), Eq.~(\ref{ii4}) becomes
\begin{eqnarray}
\Phi(w^{-1},n,1-b) + (-1)^{n}\,w\,\Phi(w,n,b)& = &   \nonumber  \\
 & & \hspace{-140pt} \frac{\pi(-1)^{n-1}}{(n-1)!}\,\frac{\partial^{n-1}}{\partial (1-b)^{n-1}}\left(w^{1-b}\left(\cot(\pi-\pi b)
 -{\rm sgn}(\varphi)\,{i}\right)\right).   \label{v1}
\end{eqnarray}
Rearrangement of terms and use of the representation (\ref{i1}) for $\Phi(w^{-1},n,1-b)$ gives (\ref{ii5}). $\square$

\medskip
\noindent {\bf Remark 3.} It is trivial to check that our expansion (\ref{ii5}), whenever the arguments of $\Phi$ make it applicable,
satisfies the well known identities \cite[Eqs. (7)--(9)]{guil}
\begin{eqnarray}
\Phi(z,s,a+1) & = & \frac{1}{z}\left(\Phi(z,s,a)-\frac{1}{a^s}\right),   \label{v3}  \\
\Phi(z,s-1,a) & = & \left(a+z\frac{\partial}{\partial z}\right)\Phi(z,s,a)\,,      \label{v4}  \\
\Phi(z,s+1,a) & = & -\,\frac{1}{s}\,\frac{\partial}{\partial a}\Phi(z,s,a)\,,     \label{v5}
\end{eqnarray}
stemming from the series representation (\ref{i1}), and the partial differential equation \cite[Eq.~(1.7)]{laga}
\begin{equation}
\left(z\,\frac{\partial}{\partial z}\,\frac{\partial}{\partial a}+a\,\frac{\partial}{\partial a}+ s\right)\Phi(z,s,a)=0\,, \label{v6}
\end{equation}
which results from  composition of the identities (\ref{v4}) and (\ref{v5}).

\medskip
\noindent {\bf Remark 4.} The expansion (\ref{ii5}) is apparently not valid in the case of positive integer values of $b=N=1, 2, \ldots$, since the right-hand side shows singular terms. Nevertheless, Eq. (\ref{ii5}) becomes meaningful for $b=N$ if we interpret the right-hand side as its limit for $b\to N$, that is,
\begin{eqnarray}
\Phi(w,n,N) & = & \lim_{t\to -N}\left\{\frac{\pi}{(n-1)!}\frac{\partial^{n-1}}{\partial t^{n-1}}
\left(-w^{t}\cot(\pi t)\right)-\frac{w^{-N}}{(-t-N)^n}\right\}  \nonumber  \\
& & \hspace{-40pt}+\,w^{-N}\left({\rm sgn}(\varphi)\,\frac{{i}\,\pi\,(\ln w)^{n-1}}{(n-1)!} -\sum_{k=1}^{N-1}\frac{w^k}{k^n}-
(-1)^n\,{\rm Li}_n(w^{-1})\right),    \label{vi5}
\end{eqnarray}
or, equivalently,
\begin{eqnarray}
\Phi(w,n,N) & = & w^{-N}\Bigg( \lim_{\varepsilon\to 0}\left\{\frac{\pi}{(n-1)!}\frac{\partial^{n-1}}{\partial \varepsilon^{n-1}}
\big(-w^{\varepsilon}\cot(\pi \varepsilon)\big)-\frac{(-1)^n}{(\varepsilon)^n}\right\} \nonumber  \\
& &  +\,{\rm sgn}(\varphi)\,\frac{{i}\,\pi\,(\ln w)^{n-1}}{(n-1)!} -\sum_{k=1}^{N-1}\frac{w^k}{k^n}-(-1)^n\,{\rm Li}_n(w^{-1})\Bigg),
\label{vi6}
\end{eqnarray}
where ${\rm Li}_n$ represents, as usual, the polylogarithm of order $n$. One obtains in this way, for the lowest values of $n$,
\begin{eqnarray}
\Phi(w,1,N) & = & w^{-N}\left( {\rm sgn}(\varphi)\,{i}\,\pi-\ln w-\sum_{k=1}^{N-1}\frac{w^k}{k}+{\rm Li}_1(w^{-1})\right),    \\
\Phi(w,2,N) & = & w^{-N}\Bigg( \frac{\pi^2}{3} + {\rm sgn}(\varphi)\,{i}\,\pi\,\ln w - \frac{1}{2}\,(\ln w)^2  \nonumber  \\
& & \hspace{130pt} -\,\sum_{k=1}^{N-1}\frac{w^k}{k^2}-{\rm Li}_2(w^{-1})\Bigg),    \\
\Phi(w,3,N) & = & w^{-N}\Bigg( \frac{\pi^2}{3}\,\ln w + {\rm sgn}(\varphi)\,\frac{{i}\,\pi}{2}\,(\ln w)^2 - \frac{1}{6}\,(\ln w)^3   \nonumber \\
& & \hspace{130pt} -\,\sum_{k=1}^{N-1}\frac{w^k}{k^3}+{\rm Li}_3(w^{-1})\Bigg),   \\
\Phi(w,4,N) & = & w^{-N}\Bigg( \frac{\pi^4}{45} + \frac{\pi^2}{6}\,(\ln w)^2 + {\rm sgn}(\varphi)\,\frac{{i}\,\pi}{6}\,(\ln w)^3 -
\frac{1}{24}\,(\ln w)^4   \nonumber   \\
& & \hspace{130pt} -\,\sum_{k=1}^{N-1}\frac{w^k}{k^4}-{\rm Li}_4(w^{-1})\Bigg),    \\
\Phi(w,5,N) & = & w^{-N}\Bigg( \frac{\pi^4}{45}\,\ln w + \frac{\pi^2}{18}\,(\ln w)^3 + {\rm sgn}(\varphi)\,\frac{{i}\,\pi}{24}\,(\ln w)^4 -
\frac{1}{120}\,(\ln w)^5  \nonumber  \\
 & & \hspace{130pt} -\,\sum_{k=1}^{N-1}\frac{w^k}{k^5}+{\rm Li}_5(w^{-1})\Bigg),
\end{eqnarray}
understanding that the sum in $k$ is void if $N=1$.

\medskip\noindent {\bf Remark 5.} Lagarias and Li \cite{laga} have coined the term ``extended polylogarithm" to refer to the function
\begin{equation}
{\rm Li}_s(z,a):=z\,\Phi(z,s,a)\,,   \label{vi7}
\end{equation}
which interpolates all polylogarithms via the parameter $s$. In the case, considered in this paper, of positive integer $s=n$, ${\rm Li}_n(z,a)$
gives a deformation of the polylogarithm  ${\rm Li}_n(z)$ with deformation parameter $a$. These authors have proved \cite[Theorem 8.1]{laga}
that ${\rm Li}_n(z,a)$, considered as a function of $z$, obeys a linear ordinary differential equation of order $n+1$ with polynomial coefficients.
Following a similar procedure, one can prove that $\Phi(z,n,a)$ is a solution of the differential equation
\begin{equation}
z\,\frac{{d}}{{d}z}(1-z)\left( z\,\frac{{d}}{{d}z}+a\right)^n\,f(z)=0\,.   \label{vi8}
\end{equation}
Linearly independent solutions of this equation are the functions of $z$
\begin{equation}
z^{-a}(\ln z)^{n-1}, \quad z^{-a}(\ln z)^{n-2}, \quad \ldots, \quad z^{-a},\quad \sum_{m=1}^\infty z^{-m}(a-m)^{-n},  \label{vi9}
\end{equation}
the last of them only for $|z|>1$. These functions constitute a basis in the space of solutions of (\ref{vi8}). Our expansion (\ref{ii5})
merely expresses, {\em mutatis mutandis}, the solution $\Phi(z,n,a)$ in such a basis.

\section{Final discussion}

From a computational point of view, the representation of the Lerch's transcendent given in Theorem 1, cannot compete against the series expansion (\ref{i1}). Nevertheless, written in the form (\ref{iii11}), it becomes very useful to compute principal value integrals appearing in different branches of Physics, especially in high energy scattering of elementary particles. More interesting, in what concerns Number Theory, is the property unveiled in our Theorem 2, which has its counterparts in analogous properties for the Hurwitz zeta function and the gamma function. Our third result, namely Corollary 1, provides an expansion of the Lerch's transcendent which, due to its algebraic simplicity and computational efficiency, may be favourably compared with the known expansions for large values of the first argument. It has, however, the drawback of being restricted to positive integer values of the second argument.
Asymptotic expansions of $\Phi(z,s,a)$, for large complex $z$, large or small complex $a$, and complex $s$, have been obtained by
Ferreira and L\'opez \cite{ferr}. On the other hand, Navas, Ruiz and Varona \cite{nava} have studied $\Phi(z,s,a)$ as a function of the complex variable $a$, with complex parameters $z$ and $s$, and shown its asymptotic behaviour  for $\Re s\to -\infty$ and $|\Im s|$ bounded. In our case, the restriction of $s$, in $\Phi(z,s,a)$, to positive integer values $s=n=1, 2, \ldots ,$
has been necessary to obtain the representation (\ref{ii3}). The question arises if this representation and its consequences are valid also for more general values of $s$. A conclusive answer to this issue would require methods of fractional calculus which lie out of the scope of this paper.

\section*{Acknowledgements}

This work has been supported by Conselho Nacional de Desenvolvimento Cient\'{\i}fico e Tecnol\'ogico  (CNPq, Brazil, grant 30445/2004-5)
and by Departamento de Ciencia, Tecnolog\'{\i}a y Universidad del Gobierno de Arag\'on (Project E24/1) and Ministerio de
Ciencia e Innovaci\'on (Project MTM2009-11154).

\end{document}